\documentclass[12pt]{article}

\usepackage{amsmath,amsthm,amssymb}

\theoremstyle{plain} 
\newtheorem{thm}{Theorem}

\newtheorem{lem}[thm]{Lemma}

\theoremstyle{definition}
\newtheorem{defn}{Definition}

\title{A correct proof of the McMorris-Powers' theorem on the consensus 
of phylogenies}
\author{Bhalchandra D. Thatte \\
\small Allan Wilson Centre for Molecular Ecology and Evolution \\
\small Massey University, Palmerston North, New Zealand \\
\small \texttt{b.thatte@massey.ac.nz}
}
\date{}

\begin{document}
\maketitle{}
\begin{abstract}
McMorris and Powers proved an Arrow-type theorem on phylogenies
given as collections of quartets. There is an error in one of the
main lemmas used to prove this theorem. However, this lemma
(and thereby the theorem) is still true, and we provide a
corrected proof.
\end{abstract}
\section {Theorem of McMorris and Powers}
\label{intro}
In 1951, K. Arrow \cite{arrow1951} proved his impossibility theorem
for the aggregation of weak orders. In the
past thirty years, there has been increasing interest in
studying limitations and possibilities in
aggregation or consensus models in many areas of science
other than the social choice theory. One such study
that is of interest to phylogeneticists is an
impossibility theorem about trees by McMorris and Powers,
(see \cite{mp1993}). The purpose of this note is to point out
an error in the original proof, and propose
a workaround.

The notation and definitions summarised below closely
follow Day and McMorris \cite{dm2003}. Let $S$ be an $n$-element set,
where $n \geq 5$. A {\em phylogeny} on $S$ is an unrooted
tree with no vertex of degree 2, and exactly $n$ vertices
of degree 1 ({\em leaves}), each labelled by a distinct
element of $S$. Let $\mathcal{P}$ denote the set of all
phylogenies on $S$. Let $w,x,y,z \in S$. We say that
the configuration ({\em quartet}) $wx|yz$ is in phylogeny
$T$ if the path from $w$ to $x$ has no vertex in common with
the path from $y$ to $z$. If the $w-x$ and $y-z$ paths have
exactly one vertex in common then we say that the quartet
$wxyz$ is in $T$. Any four elements
$w,x,y,z$ occur in $T$ in one of the four configurations
$wx|yz$, $wy|zx$, $wz|xy$ (called the {\em resolved
quartets}) and $wxyz$ (called an {\em unresolved quartet}).
Since a tree is uniquely determined by its collection
of quartets, (see \cite{cs1981}), we overload the
notation $T$ to sometimes denote the set of quartets
$q(T)$ of $T$.
For tree $T \in \mathcal{P}$ and $X\subseteq S$, $T|_X$
(the {\em restriction } of $T$ to $X$) denotes the set
of quartets of $T$ made up entirely with elements from
$X$. Similarly the restriction of a profile 
$P = (T_1, T_2, \ldots , T_k)\in \mathcal{P}^k$
to $X$ is simply $P|_X = (T_1|_X, T_2|_X, \ldots , T_k|_X)$.
For $w,x,y,z \in S$, the set of individuals $K$, $|K| = k$,
and a consensus rule $C:\mathcal{P}^k\rightarrow \mathcal{P}$,
the following shortcut notation is used:
$wx|yz \in q(T_i)$ is denoted by $wxT_iyz$, 
$wxyz \in q(T_i)$ is denoted by $wxyzT_i$,
$(\forall i \in I\subseteq K)(wxT_iyz)$ is denoted by $wxT_Iyz$,
and $(\forall i \in I\subseteq K)(wxyzT_i)$ is denoted by $wxyzT_I$.

\begin{defn} Let $C:\mathcal{P}^k\rightarrow \mathcal{P}$ be a
consensus function. The notions of Dictatorship (Dct),
Independence (Ind) and Pareto Optimality (PO) are defined as:\\
$ {\bf Dct:} (\exists j \in K)(\forall w,x,y,z \in S)
(\forall P\in \mathcal{P}^k)
(wxT_jyz \Rightarrow wxC(P)yz) $ \\
$ {\bf Ind:} (\forall X\subseteq S)(\forall P, P'\in \mathcal{P}^k)
(P|_X = P'|_X \Rightarrow C(P)|_X = C(P')|_X) $\\
$ {\bf PO:} (\forall w,x,y,z \in S)(\forall P \in \mathcal{P}^k)
(wxT_Kyz \Rightarrow wxC(P)yz)$\\
\end{defn}

In \cite{mp1993}, McMorris and Powers proved the following Arrow-type
Theorem.

\begin{thm}
\label{thm-arrow}
Let $C$ be a consensus rule on $\mathcal{P}$. $C$ satisfies
Dct iff it satisfies Ind and PO.
\end{thm}

Proof of this theorem follows Sen's strategy,
(see \cite{sen1970} or \cite{dm2003}),
of defining an \mbox{appropriate} notion of {\em decisiveness}
for a group of individuals, (which, informally speaking,
says, if a group of individuals has a certain preference
then the consensus profile also imposes the same preference),
followed by an {\em invariance } lemma for decisiveness
(which, informally speaking, says, if a group of individuals
is decisive about one quartet then the group is decisive
about all quartets).
One then refines the decisive set to prove that there
exists a singleton decisive set, which is a dictatorial
situation.

In the following, definitions of decisiveness are
followed by the proof of the invariance lemma as
presented by Day and McMorris in \cite{dm2003}.
The proof in \cite{mp1993} is based on the same
argument. An error in their proof is then discussed.
A new proof is presented in the next section.

\begin{defn}
\label{def-almost}
Let $C:\mathcal{P}^k\rightarrow \mathcal{P}$ be a consensus
rule, $I\subseteq K$, and $wx|yz$ a quartet.
$I$ is called almost decisive for $wx|yz$, denoted by $\hat{U}^I_{wx|yz}$,
if \\
$(\forall P\in \mathcal{P}^k)
(wxT_Iyz \wedge wxyzT_{K\backslash I}\Rightarrow wxC(P)yz)$.
$I$ is called almost decisive if it is almost decisive for
all resolved quartets.
The family of almost decisive subsets of $K$ is denoted by $\hat{U}_C$.
$I$ is called decisive for $wx|yz$, (denoted by $U^I_{wx|yz}$)
if $(\forall P\in \mathcal{P}^k)(wxT_Iyz \Rightarrow wxC(P)yz)$.
$I$ is called decisive if it is decisive for all resolved quartets.
The family of decisive subsets of $K$ is denoted by $U_C$.
\end{defn}

\begin{lem} {\bf (lemma 3.34 in \cite{dm2003})}
\label{lem-mp-invar}
Let $C:\mathcal{P}^k\rightarrow \mathcal{P}$ be a consensus
rule that satisfies Ind and PO, and $I \subseteq K$. Then
\begin{eqnarray} 
(\exists a,b,c,d \in S)(\hat{U}^I_{ab|cd})\Longleftrightarrow I\in \hat{U}_C 
\end{eqnarray}
\begin{eqnarray}
(\exists a,b,c,d \in S)(U^I_{ab|cd})\Longleftrightarrow I\in U_C
\end{eqnarray}
\end{lem}

\noindent {\bf Proof} Proof presented here is
based on \cite{dm2003}. Proof of (1) in \cite{dm2003} is
correct, but is presented here for a later reference.\\
{\bf Proof of (1)} Let $(\exists a,b,c,d \in S)(\hat{U}^I_{ab|cd})$.
Since $|S|\geq 5$, let $v\in S$ be such that
$v\not \in X=\{a,b,c,d\}$. We will show that $\hat{U}^I_{bv|cd}$.
Construct $P \in \mathcal{P}^k$ such that \\
$\{ab|cd,ab|cv,ab|dv,av|cd,bv|cd\}\subseteq T_I$, and
$\{abcd, av|bc,av|bd,av|cd,bcdv\}\subseteq T_{K\backslash I}$.
$P$ is otherwise unconstrained. Since $\hat{U}^I_{ab|cd}$,
$ab|cd\in C(P)$. By PO, $av|cd \in C(P)$. Therefore,
$bv|cd \in C(P)$. Therefore, by Ind, $\hat{U}^I_{bv|cd}$. By
\mbox{trivial} variants of this argument, $\hat{U}^I_{wx|yz}$ is
obtained for each $wx|yz$ other than $ab|cd$, thus proving
$I\in \hat{U}_C$. The converse is trivial.$\Box$\\
{\bf Proof of (2) (original proof)}
Let $(\exists a,b,c,d \in S)(U^I_{ab|cd})$.
Then $I$ is also almost decisive for $ab|cd$, so by (1)
it is almost decisive for all resolved quartets.
We would like to prove for any $P\in \mathcal{P}^k$ and
$X = \{w,x,y,z\} \subseteq S$, that $wxT_Iyz \Rightarrow wxC(P)yz$.
Since $|S| \geq 5$, we can select $v\not \in X$, and
construct a profile $P'$ such that $\{wx|yz,wx|vy, wx|vz, vwyz,vxyz\}
\subseteq T'_I$, and $\{vwxy, vwxz\}\subseteq T'_{K\backslash I}$,
and $P|_X = P'|_X$. $P'$ is otherwise unconstrained.
From (1), we have $\{wx|vy, wx|vz\}\subseteq C(P')$.
Therefore, $wx|yz \in C(P')$. But $P|_X = P'|_X$,
so $wx|yz \in C(P)$ as required. The converse is trivial. $\Box$

There is an error in the nontrivial direction of the proof of (2)
above. Profile $P'$ is chosen such that
$\{vwxy, vwxz\}\subseteq T'_{K\backslash I}$
and $P|_X = P'|_X$. This implies
$\{wxyz, wx|yz\} \cap T_{K\backslash I}\neq \emptyset$.
This means, if $P$ is such that $wy|xz\in T_{K\backslash I}$
or  $wz|xy\in T_{K\backslash I}$ then no choice of
$P'$ such that $P|_X = P'|_X$, can meet the
requirement $\{vwxy, vwxz\}\subseteq T'_{K\backslash I}$
of the construction.
Although the result of the lemma is correct,
a complete proof requires more complex arguments than
the ones provided by McMorris and Powers in their original
proof. In the next section, complete arguments will be
provided.

\section {Invariance Lemmas}
The proof presented below requires four different levels of
decisiveness, the first one being equivalent to almost
decisiveness defined above, and the last one being
the decisiveness defined above. Most proofs below follow
the line of argument that can be summarised thus: we have
a profile $P$ containing a certain configuration
on $X\subseteq S$. We want to prove that the configuration
also occurs in the consensus profile $C(P)$.
We construct a profile $P'$ that agrees
with $P$ when restricted to $X$.
Moreover, $P'$ allows us to resolve the configuration in
$C(P')$ using weaker notions of decisiveness
and their invariance lemmas. Then by Ind, we claim that $C(P)$
contains the configuration.

\begin{defn} {\bf (type-A decisiveness)}
Let $C:\mathcal{P}^k\rightarrow \mathcal{P}$ be a consensus
rule, $I\subseteq K$, and $X=\{a,b,c,d\}\subset S$.
$I$ is called decisive-A for $ab|cd$, denoted by
$A(I:C,ab|cd)$, if
\begin{eqnarray*}
(\forall P\in \mathcal{P}^k)
(i \in I \Rightarrow ab|cd \in T_i)
(i \in K\backslash I \Rightarrow abcd\in T_i)
\\ \Rightarrow (ab|cd \in C(P))
\end{eqnarray*}
$I$ is called decisive-A, denoted by $A(I:C)$,
if it is decisive-A for all resolved quartets.
\end{defn}

\begin{lem}
\label{lem-A}
Let $C$ be a consensus rule satisfying Ind and PO,
and $I\subseteq K$.
\begin{displaymath}
(\exists a,b,c,d \in S)(A(I:C,ab|cd)) \Rightarrow A(I:C)
\end{displaymath}
\end{lem}
\noindent {\bf Proof} The notion of decisive-A sets
is equivalent to the almost decisiveness in definition
~\ref{def-almost}. Also, this lemma is equivalent to
the {\em only if} part of (1) of Lemma ~\ref{lem-mp-invar}.
So we skip the proof.

\begin{defn} {\bf (type-B decisiveness)}
Let $C:\mathcal{P}^k\rightarrow \mathcal{P}$ be a consensus
rule, $I\subseteq K$, and $X=\{a,b,c,d\}\subset S$.
$I$ is called decisive-B for $ab|cd$ if
\begin{eqnarray*}
(\forall P\in \mathcal{P}^k)
(i \in I \Rightarrow ab|cd \in T_i)
(i \in K\backslash I \Rightarrow \{abcd, ab|cd\}\cap T_i \neq \emptyset)
\\ \Rightarrow (ab|cd\in C(P))
\end{eqnarray*}
$I$ is called decisive-B, denoted by $B(I:C)$,
if it is decisive-B for all resolved quartets.
\end{defn}

\begin{lem}
\label{lem-B}
Let $C$ be a consensus rule satisfying Ind and PO,
and $I\subseteq K$. Then $A(I:C)\Rightarrow B(I:C)$.
\end{lem}

\noindent {\bf Proof} Let $I$ be a decisive-A set.
Let $w,x,y,z\in S$, and $P= (T_1, T_2, \ldots, T_k)$ be
any profile satisfying
\begin{enumerate}
\item $(i \in I \Rightarrow wx|yz \in T_i)$,
\item $(i \in K\backslash I \Rightarrow 
\{wxyz, wx|yz\}\cap T_i \neq \emptyset)$.
\end{enumerate}
Since $|S| \geq 5$, we can select $v\not \in X$, and
construct a profile $P'$ such that 
\begin{enumerate}
\item $\{wx|yz,wx|vy, vwyz,vxyz\} \subseteq T'_I$, 
\item $\{vwxy, vwxz, wxyz\}\subseteq T'_i$ 
whenever $wxyz \in T_i$ and $i \in K\backslash I$,
\item $\{vwxy, vwxz, wx|yz\}\subseteq T'_i$ whenever
$wx|yz \in T_i$ and $i \in K\backslash I$.
\end{enumerate}
$P'$ is otherwise unconstrained. This satisfies $P|_X = P'|_X$.
If $\{wx|yz,wx|vy\}\subseteq T'_I$ then we have $wx|vz \in T'_I$. 
Since $I$ is decisive-A, $\{wx|vy, wx|vz\}\subseteq C(P')$.
Therefore, $wx|yz \in C(P')$. By $P|_X = P'|_X$, $wx|yz \in C(P)$.
Since $w,x,y,z$ are arbitrary, $I$ is decisive-B. $\Box $

This is in fact what was proved in the original proof of
part (2) of Lemma ~\ref{lem-mp-invar}.
Note that this is weaker than the full decisiveness that we desire.

\begin{defn} {\bf (type-C decisiveness)}
Let $C:\mathcal{P}^k\rightarrow \mathcal{P}$ be a consensus
rule, $I\subseteq K$, and $X=\{a,b,c,d\}\subset S$.
$I$ is called decisive-C for $ab|cd$ if
\begin{eqnarray*}
(\forall P\in \mathcal{P}^k)
(i \in I \Rightarrow ab|cd \in T_i)
(i \in K\backslash I \Rightarrow \{abcd, ab|cd, ac|bd\}\cap T_i \neq \emptyset)
\\
\Rightarrow (ab|cd\in C(P))
\end{eqnarray*}
$I$ is called decisive-C, denoted by $C(I:C)$,
if it is decisive-C for all resolved quartets.
\end{defn}
\begin{lem}
\label{lem-C}
Let $C$ be a consensus rule satisfying Ind and PO,
and $I\subseteq K$. Then $A(I:C)\Rightarrow C(I:C)$.
\end{lem}

\noindent {\bf Proof}
Let $X = \{w,x,y,z\}\subseteq S$.
Let $P = (T_1, T_2, \ldots, T_k)$ be a profile such that
$wx|yz \in T_i \forall i\in I$, and
$\{wxyz, wx|yz, wy|xz\}\cap T_i \neq \emptyset$
whenever $i \in K\backslash I$. 
We would like to prove that $wx|yz \in C(P)$.
Construct a profile $P' = (T'_1, T'_2, \ldots, T'_k)$ such that
\begin{enumerate}
\item $\{wx|yz, wx|vy, xy|vz, wy|vz\} \subseteq T'_i$
whenever $i \in I$.
\item $\{wyvx, wyvz, wv|xz, vy|xz\}\subseteq T'_i$
whenever $wy|xz \in T_i$ and $i \in K\backslash I$.
\item $\{wx|yz, wx|vy, xy|vz, wy|vz\}\subseteq T'_i$
whenever $wx|yz \in T_i$ and $i \in K\backslash I$.
\item $\{wxyz, wxyv, wyvz, wxvz, xyzv\}\subseteq T'_i $ whenever
$wxyz \in T_i$ and  $i \in K\backslash I$.
\end{enumerate}
$P'$ is otherwise unconstrained. Clearly, $P|_X = P'|_X$.
When $\{wx|yz,wx|vy\}\subseteq T'_i$, we have $wx|vz \in T'_i$.
Similarly, if $\{wv|xz, vy|xz\}\subseteq T'_i$ then  $wy|xz \in T'_i$,
and if $\{xy|vz, wy|vz\}\subseteq T'_i$ then $wx|vz  \in T'_i$.
By Lemma~\ref{lem-B}, we have $A(I:C)\Rightarrow B(I:C)$. Applying
Lemma~\ref{lem-B} to $w,y,v,z$, we have $wy|vz \in C(P')$,
and applying Lemma~\ref{lem-B} to $w,x,v,y$, we have
$wx|vy \in C(P')$. Therefore, $wx|yz \in C(P')$,
and $wx|yz \in C(P)$ by Ind. Since $w,x,y,z$ are arbitrary,
$I$ is decisive-C. $\Box $

\begin{defn} {\bf (type-D decisiveness)}
Let $C:\mathcal{P}^k\rightarrow \mathcal{P}$ be a consensus
rule, $I\subseteq K$, and $X=\{a,b,c,d\}\subset S$.
$I$ is called decisive-D (or simply decisive) for $ab|cd$
if
\begin{eqnarray*}
(\forall P\in \mathcal{P}^k)
(i \in I \Rightarrow ab|cd \in T_i) \Rightarrow (ab|cd\in C(P))
\end{eqnarray*}
$I$ is called decisive-D, or simply decisive, denoted by $D(I:C)$
if it is decisive for all resolved quartets.
\end{defn}

\begin{lem}
\label{lem-D}
Let $C$ be a consensus rule satisfying Ind and PO,
and $I\subseteq K$. Then $A(I:C)\Rightarrow D(I:C)$.
\end{lem}
\noindent {\bf Proof } Let $A(I:C)$, so by previous lemmas, we have
$B(I:C)$ and $C(I:C)$.
Let $X = \{w,x,y,z\}\subseteq S$.
Let $P = (T_1, T_2, \ldots, T_k)$ be a profile such that
$wx|yz \in T_i \forall i\in I$. P is unconstrained otherwise.
Construct a profile $P' = (T'_1, T'_2, \ldots, T'_k)$ such that
\begin{enumerate}
\item $\{wx|yz, wx|vy, wyvz, xyvz\} \subseteq T'_i$
whenever $i \in I$.
\item $\{wyvx, wyvz, wv|xz, vy|xz\}\subseteq T'_i$
whenever $wy|xz \in T_i$ and $i \in K\backslash I$.
\item $\{wzvx, wzvy, vz|xy, wv|xy\}\subseteq T'_i$
whenever $wz|xy \in T_i$ and $i \in K\backslash I$.
\item $\{wxvy, wxvz, wv|yz, xv|yz\}\subseteq T'_i $
whenever $wx|yz \in T_i$ and  $i \in K\backslash I$.
\item $\{wxvy, wxvz, wxyz, wvyz, xvyz\}\subseteq T'_i $
whenever $wxyz \in T_i$ and  $i \in K\backslash I$.
\end{enumerate}
$P'$ is otherwise unconstrained.
Clearly, $P|_X = P'|_X$.
When $\{wx|yz,wx|vy\}\subseteq T'_i$, we have $wx|vz \in T'_i$.
Similarly, when $\{wv|xz, vy|xz\}\subseteq T'_i$, we have $wy|xz \in T'_i$,
and if  $\{vz|xy, wv|xy\}\subseteq T'_i$ then $wz|xy\in T'_i$,
and if $\{wv|yz, xv|yz\}\subseteq T'_i $ then $wx|yz \in T'_i$.
By Lemma~\ref{lem-C}, we have $\{wx|vy, wx|vz\} \subseteq C(P')$,
which implies $wx|yz \in C(P')$ and,
by Ind, $wx|yz\in C(P)$. $\Box$

This lemma implies part (2) of Lemma ~\ref{lem-mp-invar}.

\section* {Acknowledgments} I would like to thank Mike
Steel for making several suggestions to improve the
presentation of this manuscript.

\end{document}